\documentclass[12pt]{article}
\usepackage{amssymb,amsmath}
\newtheorem{example}{ Example}[section]
\newtheorem{proposition}{Proposition}[section]
\newtheorem{theorem}{Theorem}[section]

\newtheorem{remark}{Remark}[section]

\numberwithin{equation}{section}

\def\eps{\varepsilon}
\begin{document}

\begin{center}{\large\sc
Multiple solutions for a Neumann system involving subquadratic
nonlinearities
}\\
\vspace{0.5cm}

{\large  Alexandru Krist\'aly\footnote{Research supported by CNCSIS
grant PCCE-55/2008 "Sisteme diferen\c tiale \^\i n analiza
neliniar\u a \c si aplica\c tii", by the J\'anos Bolyai Research
Scholarship of the Hungarian Academy of Sciences, and by Slovenian
Research Agency grants P1-0292-0101 and J1-2057-0101.}\\
  {\normalsize Department of Economics, Babe\c s-Bolyai University, Str. Teodor Mihali, nr. 58-60, 400591
Cluj-Napoca, Romania
}\\
\vspace{0.5cm}
 Du\v san Repov\v s\\  {\normalsize Faculty  of
 Mathematics and Physics, and Faculty of Education, University of Ljubljana, Jadranska 19, 1000 Ljubljana,  Slovenia
 }\\
   }
\end{center}


\begin{abstract}
\noindent {\footnotesize
In this paper we consider the model semilinear Neumann system
$$\left\{
\begin{array}{lll}
-\Delta u+a(x)u=\lambda c(x) F_u(u,v)& {\rm in} & \Omega,\\
-\Delta v+b(x)v=\lambda c(x) F_v(u,v)& {\rm in} & \Omega,\\
\frac{\partial u}{\partial \nu}=\frac{\partial v}{\partial \nu}=0 &
{\rm on} & \partial\Omega,
\end{array}\right.
 \eqno{(N_\lambda)}
 $$
 where $\Omega\subset \mathbb R^N$ is a smooth open bounded domain, $\nu$ denotes the outward unit normal to $\partial \Omega$, $\lambda\geq 0$ is a parameter, $a,b,c\in L_+^\infty(\Omega)\setminus\{0\},$ and $F\in C^1(\mathbb{R}^2,\mathbb{R})\setminus\{0\}$ is a nonnegative function which is
 subquadratic at infinity. Two nearby numbers are determined in explicit forms, $\underline \lambda$ and $\overline \lambda$ with
 $ 0<\underline\lambda\leq \overline \lambda$, such that for every $0\leq \lambda<\underline \lambda$, system $(N_\lambda)$
 has only the trivial pair of solution, while for every $\lambda>\overline \lambda$, system
 $(N_\lambda)$ has at least two distinct nonzero pairs
 of solutions.
}
\end{abstract}


\section{Introduction}
Let us consider the quasilinear Neumann system
$$\left\{
\begin{array}{lll}
-\Delta_p u+a(x)|u|^{p-2}u=\lambda c(x) F_u(u,v)& {\rm in} & \Omega,\\
-\Delta_q v+b(x)|v|^{q-2}v=\lambda c(x) F_v(u,v)& {\rm in} & \Omega,\\
\frac{\partial u}{\partial \nu}=\frac{\partial v}{\partial \nu}=0 &
{\rm on} & \partial\Omega,
\end{array}\right.
 \eqno{(N^{p,q}_\lambda)}
 $$
 where $p,q>1;$ $\Omega\subset \mathbb R^N$ is a smooth open bounded domain; $\nu$ denotes the outward unit normal to $\partial \Omega$;
 $a,b,c\in L^\infty(\Omega)$ are some functions;  $\lambda\geq 0$ is a
 parameter; and $F_u$ and $F_v$ denote the partial derivatives of $F\in C^1(\mathbb{R}^2,\mathbb{R})$ with
 respect to the first and second variables, respectively.

Recently, problem $(N^{p,q}_\lambda)$ has been considered by several
authors. For instance, under suitable assumptions on $a,b,c$ and
$F$,  El Manouni and Kbiri Alaoui \cite{EM-KA} proved the existence
of an interval $A\subset (0,\infty)$ such that $(N^{p,q}_\lambda)$
has at least three solutions whenever $\lambda\in A$ and $p,q>N$.
Lisei and Varga \cite{LV} also established the existence of at least
three solutions for the system $(N^{p,q}_\lambda)$ with
nonhomogeneous and nonsmooth Neumann boundary conditions. Di Falco
\cite{DiFalco} proved the existence of infinitely many solutions for
$(N^{p,q}_\lambda)$ when the nonlinear function $F$ has a suitable
oscillatory behavior. Systems similar to $(N^{p,q}_\lambda)$ with
the Dirichlet boundary conditions were also considered by Afrouzi
and Heidarkhani \cite{AH-1, AH-2}, Boccardo and de Figueiredo
\cite{deFig-Boccardo}, Heidarkhani and Tian \cite{HT}, Li and Tang
\cite{LT}, see also references therein.

The aim of the present paper is to describe a new phenomenon for
Neumann systems when the nonlinear term has a subquadratic growth.
In order to avoid technicalities, instead of the quasilinear system
$(N^{p,q}_\lambda)$, we shall consider the semilinear problem
$$\left\{
\begin{array}{lll}
-\Delta u+a(x)u=\lambda c(x) F_u(u,v)& {\rm in} & \Omega,\\
-\Delta v+b(x)v=\lambda c(x) F_v(u,v)& {\rm in} & \Omega,\\
\frac{\partial u}{\partial n}=\frac{\partial v}{\partial n}=0 & {\rm
on} & \partial\Omega.
\end{array}\right.
 \eqno{(N_\lambda)}
 $$

We assume that the nonlinear term $F\in
C^1(\mathbb{R}^2,\mathbb{R})$ satisfies the following properties:
\begin{itemize}
  \item[{ $(\bf F_+)$}] $F(s,t)\geq 0$ for every $(s,t)\in \mathbb
  R^2$, $F(0,0)=0,$ and $F\not\equiv 0;$
  \item[{ $(\bf F_0)$}] 
  $\lim_{(s,t)\to (0,0)}\frac{F_s(s,t)}{{|s|}+{|t|}}=\lim_{(s,t)\to (0,0)}\frac{F_t(s,t)}{{|s|}+{|t|}}=0;$
  \item[{ $(\bf F_\infty)$}] 
   $\lim_{|s|+|t|\to \infty}\frac{F_s(s,t)}{{|s|}+{|t|}}=\lim_{|s|+|t|\to \infty}\frac{F_t(s,t)}{{|s|}+{|t|}}=0.$
\end{itemize}

\begin{example}\rm\label{pelda}
A typical nonlinearity which fulfils hypotheses $(\bf F_+)$, $(\bf
F_0)$ and $(\bf F_\infty)$ is $F(s,t)=\ln(1+s^2t^2).$
\end{example}
We also introduce the set $$\Pi_+(\Omega)=\{a\in L^\infty(\Omega):
{\rm essinf}_\Omega a>0\}.$$ For $a,b,c\in \Pi_+(\Omega)$ and for
$F\in C^1(\mathbb R^2,R)$ which fulfils the hypotheses $(\bf F_+)$,
$(\bf F_0)$ and $(\bf F_\infty)$, we define  the numbers
$$s_F=2\|c\|_{L^1}\max_{(s,t)\neq
(0,0)}\frac{F(s,t)}{\|a\|_{L^1}{s^2}+\|b\|_{L^1}t^2},\ {\rm and}\
S_F=\max_{(s,t)\neq
(0,0)}\frac{|sF_s(s,t)+tF_t(s,t)|}{{\|c/a\|_{L^\infty}^{-1}s^2}+{\|c/b\|_{L^\infty}^{-1}t^2}}.$$
Note that these numbers are finite, positive and $S_F\geq s_F$, see
Proposition \ref{elso-prop} (here and in the sequel,
$\|\cdot\|_{L^p}$ denotes the usual norm of the Lebesgue space
$L^p(\Omega)$, $p\in [1,\infty]$). Our main result reads as follows.

\begin{theorem}\label{fotetel-NA} Let  $F\in C^1(\mathbb{R}^2,\mathbb{R})$ be
a function which satisfies $(\bf F_+)$, $(\bf F_0)$ and $(\bf
F_\infty)$, and $a,b,c\in \Pi_+(\Omega)$. Then, the following
statements hold.
\begin{enumerate}
\item[{\rm (i)}] For every
$0\leq\lambda<S_F^{-1}$, system $(N_\lambda)$ has only the trivial
pair of solution.
\item[{\rm (ii)}] For every $\lambda>s_F^{-1}$, system $(N_\lambda)$ has at least two
distinct, nontrivial pairs of solutions
$(u_\lambda^i,v_\lambda^i)\in H^1(\Omega)^2$, $i\in \{1,2\}$.
 \end{enumerate}
\end{theorem}

\begin{remark}\rm (a) A natural question arises which is still open: how many solutions exist for $(N_\lambda)$  when
$\lambda\in [S_F^{-1},s_F^{-1}]$? Numerical experiments show that
$s_F$ and $S_F$ are usually not far from each other, although their
origins are independent. For instance, if $a=b=c$, and $F$ is from
Example \ref{pelda}, we have $s_F\approx 0.8046$ and $S_F=1.$

(b) Assumptions $(\bf F_+)$, $(\bf F_0)$ and $(\bf F_\infty)$ imply
that there exists $c>0$ such that
\begin{equation}\label{F-re-vonatkozo-szukvadratikus}
    0\leq F(s,t)\leq c(s^2+t^2) \ {\rm for\ all}\ (s,t)\in \mathbb R^2,
\end{equation}
i.e., $F$ has a subquadratic growth. Consequently, Theorem
\ref{fotetel-NA} completes the results of several papers where $F$
fulfils the Ambrosetti-Rabinowitz condition, i.e., there exist
$\theta>2$ and $r>0$ such that
\begin{equation}\label{AR}
0<\theta F(s,t)\leq s F_s(s,t)+t F_t(s,t)\ {\rm for\ all}\
|s|,|t|\geq r.
\end{equation}
Indeed, (\ref{AR}) implies that for some $C_1,C_2>0$, one has
$F(s,t)\geq C_1(|s|^\theta+|t|^\theta)$ for all $|s|,|t|> C_2.$
\end{remark}

The next section contains some auxiliary notions and results, while
in Section \ref{sect-bizonyitas} we prove Theorem \ref{fotetel-NA}.
First, a direct calculation proves (i), while a very recent three
critical points result of  Ricceri \cite{Ri-1} provides the proof of
(ii).

\section{Preliminaries}

A solution for $(N_\lambda)$ is a pair $(u,v)\in H^1(\Omega)^2$ such
that
\begin{equation}\label{va}
    \left\{
\begin{array}{lll}
\int_\Omega (\nabla u \nabla \phi+a(x)u\phi)dx=\lambda \int_\Omega c(x) F_u(u,v)\phi dx& {\rm for\ all} & \phi\in H^1(\Omega),\\
\int_\Omega (\nabla v \nabla \psi+b(x)v\psi)dx=\lambda \int_\Omega
c(x) F_v(u,v)\psi dx& {\rm for\ all} & \psi\in H^1(\Omega).
\end{array}\right.
\end{equation}

Let $a,b,c\in \Pi_+(\Omega).$ We associate to the system
$(N_\lambda)$ the energy functional $I_\lambda:H^1(\Omega)^2\to
\mathbb R$ defined by
$$I_\lambda(u,v)=\frac{1}{2}(\|u\|_a^2+\|v\|_b^2)-\lambda \mathcal F(u,v),$$
where $$\|u\|_a=\left(\int_\Omega |\nabla
u|^2+a(x)u^2\right)^{1/2};\ \|v\|_b=\left(\int_\Omega |\nabla
v|^2+b(x)v^2\right)^{1/2},$$ and $$\mathcal F(u,v)=\int_\Omega
c(x)F(u,v).$$ It is clear that $\|\cdot \|_a$ and $\|\cdot \|_b$ are
equivalent to the usual norm on $H^1(\Omega).$ Note that if $F\in
C^1(\mathbb R^2,R)$  verifies the hypotheses $(\bf F_0)$ and $(\bf
F_\infty)$ (see also relation
(\ref{F-re-vonatkozo-szukvadratikus})), the functional $I_\lambda$
is well-defined, of class $C^1$ on $H^1(\Omega)^2$ and its critical
points are exactly the solutions for $(N_\lambda)$. Since
$F_s(0,0)=F_t(0,0)=0$ from $(\bf F_0)$, $(0,0)$ is a solution of
$(N_\lambda)$ for every $\lambda\geq 0.$

In order to prove Theorem \ref{fotetel-NA} (ii), we must find
critical points for $I_\lambda.$ In order to do this, we recall the
following Ricceri-type three critical point theorem. First, we need
the following notion: if $X$ is a Banach space, we denote by
$\mathcal W_X$ the class of those functionals $E:X\to \mathbb{R}$
that possess the property that if $\{u_n\}$ is a sequence in $X$
converging weakly to $u\in X$ and $\liminf_n E(u_n)\leq E(u)$ then
$\{u_n\}$ has a subsequence strongly converging to $u$.

\begin{theorem} {\rm \cite[Theorem 2]{Ri-1}}\label{bonanno-tetel}
Let $X$ be a separable and reflexive real Banach space, let
$E_1:X\to \mathbb{R}$ be a coercive, sequentially weakly lower
semicontinuous $C^1$ functional belonging to $\mathcal W_X$, bounded
on each bounded subset of $X$ and whose derivative admits a
continuous inverse on $X^*$; and $E_2:X\to \mathbb{R}$ a $C^1$
functional with a compact derivative. Assume that $E_1$ has a strict
local minimum $u_0$ with $E_1(u_0)=E_2(u_0)=0$. Setting the numbers
\begin{equation}\label{tau}
\tau=\max\left\{0,\limsup_{\|u\|\to
\infty}\frac{E_2(u)}{E_1(u)},\limsup_{u\to
u_0}\frac{E_2(u)}{E_1(u)}\right\},
\end{equation}
\begin{equation}\label{chi}
    \chi=\sup_{E_1(u)>0}\frac{E_2(u)}{E_1(u)},
\end{equation}
assume that $\tau<\chi.$

Then, for each compact interval $[a,b]\subset (1/\chi,1/\tau)$
$($with the
 conventions
$1/0=\infty$ and $1/\infty=0$$)$ there exists  $\kappa>0$ with the
following property: for every $\lambda\in [a,b]$ and every $C^1$
functional $E_3:X\to \mathbb{R}$ with a compact derivative, there
exists $\delta>0$ such that for each $\theta\in [0,\delta],$ the
equation $$E_1'(u)-\lambda E_2'(u)-\theta E_3'(u)=0$$ admits at
least three solutions in $X$ having norm less than $\kappa.$
\end{theorem}

\noindent We conclude this section with an  observation which
involves the constants $s_F$ and $S_F$.

\begin{proposition}\label{elso-prop} Let  $F\in C^1(\mathbb{R}^2,\mathbb{R})$ be
a function which satisfies $(\bf F_+)$, $(\bf F_0)$ and $(\bf
F_\infty)$, and $a,b,c\in \Pi_+(\Omega)$. Then the numbers $s_F$ and
$S_F$ are finite, positive and $S_F\geq s_F.$
\end{proposition}

{\it Proof.} It follows by $(\bf F_0)$ and $(\bf F_\infty)$ and by
the continuity of the functions $(s,t)\mapsto
\frac{F_s(s,t)}{|s|+|t|}$, $(s,t)\mapsto \frac{F_t(s,t)}{|s|+|t|}$
away from $(0,0)$,  that there exists $M>0$ such that
$$|F_s(s,t)|\leq M(|s|+|t|)\ {\rm and}\ |F_t(s,t)|\leq M(|s|+|t|)\ {\rm for\ all}\ (s,t)\in \mathbb
R^2.$$ Consequently, a standard mean value theorem together with
$(\bf F_+)$ implies that
\begin{equation}\label{M-es-becsles}
    0\leq F(s,t)\leq 2M(s^2+t^2)\ {\rm for\ all}\ (s,t)\in \mathbb R^2.
\end{equation}
 We now prove that
\begin{equation}\label{zeroban-F}
    \lim_{(s,t)\to (0,0)}\frac{F(s,t)}{{s^2}+{t^2}}=0\ {\rm and}\ \lim_{|s|+|t|\to \infty}\frac{F(s,t)}{{s^2}+{t^2}}=0.
\end{equation}

By $(\bf F_0)$ and $(\bf F_\infty)$, for every $\eps>0$ there exists
$\delta_\eps\in (0,1)$ such that for every $(s,t)\in \mathbb R^2$
with $|s|+|t|\in (0,\delta_\eps)\cup (\delta_\eps^{-1},\infty)$, one
has
\begin{equation}\label{eps-F}
    \frac{|F_s(s,t)|}{{|s|}+{|t|}}<\frac{\eps}{4}\ {\rm and}\
\frac{|F_t(s,t)|}{{|s|}+{|t|}}<\frac{\eps}{4}.
\end{equation}
By (\ref{eps-F}) and the mean value theorem, for every $(s,t)\in
\mathbb R^2$ with $|s|+|t|\in (0,\delta_\eps)$,
 we have
\begin{eqnarray*}
  F(s,t) &=& F(s,t)-F(0,t)+F(0,t)-F(0,0) \\
   &\leq& \frac{\eps}{2}(s^2+t^2)
\end{eqnarray*}
which gives the first limit in (\ref{zeroban-F}). Now, for every
$(s,t)\in \mathbb R^2$ with
$|s|+|t|>\delta_\eps^{-1}\max\{1,\sqrt{8M/\eps}\}$, by using
(\ref{M-es-becsles}) and (\ref{eps-F}),
 we have
\begin{eqnarray*}
  F(s,t) &=& F(s,t)-
  F\left(\frac{\delta_\eps^{-1}}{|s|+|t|}s,t\right)
  +F\left(\frac{\delta_\eps^{-1}}{|s|+|t|}s,t\right)
  -F\left(\frac{\delta_\eps^{-1}}{|s|+|t|}s,\frac{\delta_\eps^{-1}}{|s|+|t|}t\right)\\ &&+F\left(\frac{\delta_\eps^{-1}}{|s|+|t|}s,\frac{\delta_\eps^{-1}}{|s|+|t|}t\right)  \\
   &\leq& \frac{\eps}{4}(|s|+|t|)^2+2M \delta_\eps^{-2}\\&\leq&
   \eps(s^2+t^2),
\end{eqnarray*}
which leads us to the second limit in (\ref{zeroban-F}).

The facts above show that the numbers $s_F$ and $S_F$ are finite.
Moreover, $s_F>0$. We now prove that $S_F\geq s_F.$ To do this, let
$(s_0,t_0)\in \mathbb R^2\setminus \{(0,0)\}$ be a maximum point of
the function
$(s,t)\mapsto\frac{F(s,t)}{\|a\|_{L^1}{s^2}+\|b\|_{L^1}t^2}.$ In
particular, its partial derivatives vanishes at $(s_0,t_0)$,
yielding
$$F_s(s_0,t_0)(\|a\|_{L^1}{s_0^2}+\|b\|_{L^1}t_0^2)=2\|a\|_{L^1}s_0F(s_0,t_0);$$
$$F_t(s_0,t_0)(\|a\|_{L^1}{s_0^2}+\|b\|_{L^1}t_0^2)=2\|b\|_{L^1}t_0F(s_0,t_0).$$
From the two relations above we obtain that
$$s_0 F_s(s_0,t_0)+t_0 F_t(s_0,t_0)=2F(s_0,t_0).$$
On the other hand, since $a,b,c\in \Pi_+(\Omega)$, we have that
$$\|c\|_{L^1}=\int_\Omega c(x)dx=\int_\Omega \frac{c(x)}{a(x)}
a(x)dx\leq \left\|\frac{c}{a}\right\|_{L^\infty}\int_\Omega
a(x)dx=\left\|\frac{c}{a}\right\|_{L^\infty}\|a\|_{L^1},$$ thus
$\|c/a\|_{L^\infty}^{-1}\leq \|a\|_{L^1}/\|c\|_{L^1}$  and in a
similar way $\|c/b\|_{L^\infty}^{-1}\leq \|b\|_{L^1}/\|c\|_{L^1}$.
Combining these inequalities with the above argument,
we conclude that $S_F\geq s_F.$  \hfill $\square$\\


\vspace{-0.5cm}
\section{Proof of Theorem \ref{fotetel-NA}}\label{sect-bizonyitas}

In this section we assume that the assumptions of Theorem
\ref{fotetel-NA} are fulfilled.\\

\noindent {\bf Proof of Theorem \ref{fotetel-NA} (i).} Let $(u,v)\in
H^1(\Omega)^2$ be a solution of $(N_\lambda).$ Choosing $\phi=u$ and
$\psi=v$ in (\ref{va}),  we obtain that
\begin{eqnarray*}
  \|u\|_a^2 + \|v\|_b^2&=& \int_\Omega (|\nabla
u|^2+a(x)u^2+|\nabla
v|^2+b(x)v^2) \\
   &=& \lambda \int_\Omega c(x)(F_u(u,v)u+F_v(u,v)v)
 \\
   &\leq & \lambda S_F \int_\Omega c(x)({\|c/a\|_{L^\infty}^{-1}u^2}+{\|c/b\|_{L^\infty}^{-1}v^2}) \\
   &\leq &\lambda S_F \int_\Omega (a(x)u^2+b(x)v^2)\\ &\leq &\lambda S_F
   (\|u\|_a^2 + \|v\|_b^2).
\end{eqnarray*}
Now, if $0\leq \lambda< S_F^{-1},$ we necessarily have that
$(u,v)=(0,0)$, which concludes the proof.
\newpage

\noindent {\bf Proof of Theorem \ref{fotetel-NA} (ii).} In Theorem
\ref{bonanno-tetel} we choose $X=H^1(\Omega)^2$ endowed with the
norm $\|(u,v)\|=\sqrt{\|u\|_a^2+\|v\|_b^2}$, and $E_1,
E_2:H^1(\Omega)^2\to \mathbb R$ defined by
$$E_1(u,v)=\frac{1}{2}\|(u,v)\|^2\ {\rm and}\ E_2(u,v)=\mathcal F(u,v).$$
It is clear that both $E_1$ and $E_2$ are $C^1$ functionals and
$I_\lambda=E_1-\lambda E_2.$ It is also a standard fact that $E_1$
is a coercive, sequentially weakly lower semicontinuous functional
which belongs to $\mathcal W_{H^1(\Omega)^2}$, bounded on each
bounded subset of $H^1(\Omega)^2$, and its derivative admits a
continuous inverse on $(H^1(\Omega)^2)^*.$ Moreover, $E_2$ has a
compact derivative since $H^1(\Omega)\hookrightarrow L^p(\Omega)$ is
a compact embedding for every $p\in (2,2^*)$.

Now, we prove that the functional $(u,v)\mapsto
\frac{E_2(u,v)}{E_1(u,v)}$ has similar properties as the function
$(s,t)\mapsto \frac{F(s,t)}{s^2+t^2}.$ More precisely, we shall
prove that
\begin{equation}\label{egy-ketto}
    \lim_{\|(u,v)\|\to 0}\frac{E_2(u)}{E_1(u)}=\lim_{\|(u,v)\|\to \infty}\frac{E_2(u)}{E_1(u)}=
    0.
\end{equation}
First, relation
 (\ref{zeroban-F}) implies that for every $\eps>0$ there exists
$\delta_\eps\in (0,1)$ such that for every $(s,t)\in \mathbb R^2$
with $|s|+|t|\in (0,\delta_\eps)\cup (\delta_\eps^{-1},\infty)$, one
has
\begin{equation}\label{eps-F-uj}
    0\leq \frac{F(s,t)}{{s^2}+t^2}<\frac{\eps}{4\max\{\|c/a\|_{L^\infty},\|c/b\|_{L^\infty}\}}.
\end{equation}
Fix  $p\in (2,2^*).$ Note that the continuous function $(s,t)\mapsto
\frac{F(s,t)}{|s|^p+|t|^p}$ is bounded on the set $\{(s,t)\in
\mathbb R^2:|s|+|t|\in [\delta_\eps,\delta_\eps^{-1}]\}$. Therefore,
for some $m_\eps>0,$ we have that in particular
$$0\leq F(s,t)\leq \frac{\eps}{4\max\{\|c/a\|_{L^\infty},\|c/b\|_{L^\infty}\}}(s^2+t^2) + m_\eps (|s|^p+|t|^p)\ \ {\rm for\ all}\ (s,t)\in \mathbb R^2.$$
Therefore, for each $(u,v)\in H^1(\Omega)^2,$ we get
\begin{eqnarray*}
 0\leq E_2(u,v) &=& \int_\Omega c(x)F(u,v)  \\
   &\leq& \int_\Omega c(x)\left[\frac{\eps}{4\max\{\|c/a\|_{L^\infty},\|c/b\|_{L^\infty}\}}(u^2+v^2) + {m_\eps} (|u|^p+|v|^p)\right] \\
   &\leq& \int_\Omega \left[\frac{\eps}{4}(a(x)u^2+b(x)v^2) + {m_\eps} c(x)(|u|^p+|v|^p)\right] \\
   &\leq& \frac{\eps}{4}\|(u,v)\|^2+ {m_\eps} \|c\|_{L^\infty}S_{p}^{p}(\|u\|_a^{p}+\|v\|_b^{p})\\
   &\leq& \frac{\eps}{4}\|(u,v)\|^2+ {m_\eps} \|c\|_{L^\infty}S_{p}^{p}\|(u,v)\|^p,
\end{eqnarray*}
where $S_{l}>0$ is the best constant in the inequality
$\|u\|_{L^l}\leq S_l \min \{\|u\|_a,\|u\|_b\}$ for every $u\in
H^1(\Omega)$, $l\in (1,2^*)$ (we used the fact that the function
$\alpha\mapsto (s^\alpha+t^\alpha)^\frac{1}{\alpha}$ is decreasing,
$s,t\geq 0$). Consequently, for every $(u,v)\neq (0,0),$ we obtain
\begin{eqnarray*}
  0\leq \frac{E_2(u,v)}{E_1(u,v)} &\leq& \frac{\eps}{2}+ {2m_\eps}
  \|c\|_{L^\infty}S_{p}^{p}\|(u,v)\|^{p-2}.
\end{eqnarray*}
Since $p>2$ and $\eps>0$ is arbitrarily small when $(u,v)\to 0,$ we
obtain the first limit from (\ref{egy-ketto}).

 Now, we fix $r\in (1,2).$ The continuous function $(s,t)\mapsto
\frac{F(s,t)}{|s|^r+|t|^r}$ is bounded on the set $\{(s,t)\in
\mathbb R^2:|s|+|t|\in [\delta_\eps,\delta_\eps^{-1}]\}$, where
$\delta_\eps\in (0,1)$ is from (\ref{eps-F-uj}). Combining this fact
with (\ref{eps-F-uj}), one can find a number $M_\eps>0$ such that
$$0\leq F(s,t)\leq \frac{\eps}{4\max\{\|c/a\|_{L^\infty},\|c/b\|_{L^\infty}\}}(s^2+t^2) + M_\eps (|s|^r+|t|^r)\ \ {\rm for\ all}\ (s,t)\in \mathbb R^2.$$
The H\"older inequality and a similar calculation as above show that
\begin{eqnarray*}
  0\leq E_2(u,v) &\leq &\frac{\eps}{4}\|(u,v)\|^2+ 2^{1-\frac{r}{2}}{M_\eps} \|c\|_{L^\infty}S_{r}^{r}\|(u,v)\|^{r}.
\end{eqnarray*}
For every $(u,v)\neq (0,0)$, we have that
\begin{eqnarray*}
  0\leq \frac{E_2(u,v)}{E_1(u,v)} &\leq& \frac{\eps}{2}+ {2^{2-\frac{r}{2}}M_\eps}
  \|c\|_{L^\infty}S_{r}^{r}\|(u,v)\|^{r-2}.
\end{eqnarray*}
 Due to the arbitrariness of $\eps>0$ and $r\in (1,2)$, by letting the limit $\|(u,v)\|\to
 \infty$,  we obtain the second relation from (\ref{egy-ketto}).

Note that $E_1$ has a strict global minimum $(u_0,v_0)=(0,0)$, and
$E_1(0,0)=E_2(0,0)=0.$ The definition of the number $\tau$ in
Theorem \ref{bonanno-tetel}, see (\ref{tau}), and the limits in
(\ref{egy-ketto}) imply that $\tau=0.$ Furthermore, since
$H^1(\Omega)$ contains the constant functions on $\Omega$, keeping
the notation from (\ref{chi}), we obtain
$$\chi=\sup_{E_1(u,v)>0}\frac{E_2(u,v)}{E_1(u,v)}\geq 2\|c\|_{L^1}\max_{(s,t)\neq
(0,0)}\frac{F(s,t)}{\|a\|_{L^1}{s^2}+\|b\|_{L^1}t^2}=s_F.$$
Therefore, applying Theorem \ref{bonanno-tetel} (with $E_3\equiv0$),
we obtain that in particular for every $\lambda\in
(s_F^{-1},\infty)$, the equation $I'_\lambda(u,v)\equiv
E'_1(u,v)-\lambda E_2'(u,v)=0$ admits at least three distinct pairs
of solutions in $H^1(\Omega)^2.$  Due to condition $(\bf F_0)$,
system $(N_\lambda)$ has the solution $(0,0)$. Therefore, for every
$\lambda>s_F^{-1}$, the system $(N_\lambda)$ has at least two
distinct, nontrivial pairs of solutions, which concludes the proof.

\begin{remark}\rm\label{remark-utolso}
The conclusion of Theorem \ref{bonanno-tetel} gives a much more
precise information about the Neumann system $(N_\lambda)$; namely,
one can see that $(N_\lambda)$ is stable with respect to small
perturbations. To be more precise, let us consider the perturbed
system
$$\left\{
\begin{array}{lll}
-\Delta u+a(x)u=\lambda c(x) F_u(u,v)+\mu d(x) G_u(u,v)& {\rm in} & \Omega,\\
-\Delta v+b(x)v=\lambda c(x) F_v(u,v)+\mu d(x) G_v(u,v)& {\rm in} & \Omega,\\
\frac{\partial u}{\partial n}=\frac{\partial v}{\partial n}=0 & {\rm
on} & \partial\Omega.
\end{array}\right.
 \eqno{(N_{\lambda,\mu})}
 $$
 where $\mu\in \mathbb R,$ $d\in L^\infty(\Omega)$, and $G\in C^1(\mathbb R^2,\mathbb R)$ is a function such that for some $c>0$ and $1<p<2^*-1,$
 $$\max\{|G_s(s,t)|,|G_t(s,t)|\}\leq c(1+|s|^p+|t|^p)\ {\rm for\ all}\ (s,t)\in \mathbb R^2.$$
 One can prove in a standard manner that $E_3:H^1(\Omega)^2\to \mathbb
 R$ defined by $$E_3(u,v)=\int_{\Omega} d(x)G(u,v)dx,$$
 is of class $C^1$ and it has a compact derivative. Thus, we may apply Theorem \ref{bonanno-tetel} in its generality to show  that
  for small enough values of $\mu$ system $(N_{\lambda,\mu})$ still
 has three distinct pairs of solutions.
\end{remark}

\end{document}